\documentclass[reqno]{amsart}
\theoremstyle{plain}
\usepackage{latexsym}
\usepackage{enumerate}
\usepackage{epsfig}
\usepackage[sc]{mathpazo}
\def\XXint#1#2#3{{\setbox0=\hbox{$#1{#2#3}{\int}$}
     \vcenter{\hbox{$#2#3$}}\kern-.5\wd0}}

\newcommand{\chara}{1\!\!1}

\newcommand{\Eb}{\mathbb{E}}

\newcommand{\Fb}{\mathbb{F}}

\newcommand{\Sbb}{\mathbb{S}}

\newcommand{\lt}{\left}
\newcommand{\rt}{\right}
\newcommand{\nl}{\newline}

\newcommand{\nn}{\nonumber}

\newcommand{\lm}{\lambda}

\newcommand{\qd}{\quad}

\newcommand{\ep}{\epsilon}

\newcommand{\wt}{\widetilde}

\newcommand{\N}{\mathrm{I\!N}}

\newcommand{\PPI}{\mathcal{P}}

\newcommand{\GI}{\mathcal{G}}

\newcommand{\SI}{\mathcal{S}}

\newcommand{\UI}{\mathcal{U}}
\newcommand{\UUI}{\mathcal{U}}
\newcommand{\OI}{\mathcal{O}}

\newcommand{\ti}{\tilde}
\newcommand{\la}{\langle}
\newcommand{\ra}{\rangle}

\newcommand{\R}{\mathrm {I\!R}}

\newcommand{\dia}{\diamondsuit}

\newcommand{\FF}{\mathbb{F}}

\newtheorem{a1}{Lemma}
\newtheorem{a2}[a1]{Theorem}

\newtheorem{a6}[a1]{Corollary}

\theoremstyle{remark}

\begin{document}

\title[On indecomposable sets with applications]
{On indecomposable sets with applications}
\author{Andrew Lorent}
\address{Mathematics Department\\University of Cincinnati\\2600 Clifton Ave.\\ Cincinnati OH 45221 }
\email{lorentaw@uc.edu}
\maketitle

\begin{abstract}

In this note we show the characteristic function of every indecomposable set $F$ in the plane is $BV$ equivalent to 
the characteristic function a closed set 
$\mathbb{F}$, i.e.\ $\|\chara_{F}-\chara_{\mathbb{F}}\|_{BV(\R^2)}=0$. We show by example 
this is false in dimension three and above. As a corollary to this result we show that 
for every $\ep>0$ a set of 
finite perimeter $S$ can be approximated by a closed subset $\Sbb_{\ep}$ with finitely many 
indecomposable components and with the property that 
$H^1(\partial^M \Sbb_{\ep}\backslash \partial^M S)=0$ and 
$\|\chara_{S}-\chara_{\Sbb_{\ep}}\|_{BV(\R^2)}<\ep$. We apply this corollary to give a 
short proof that locally quasiminimizing sets in the plane are $BV_l$ extension domains. 
\end{abstract}

\section{Introduction}

Sets of finite perimeter are the largest class of sets that permit a broad theory of analysis. They have wide 
application in the Calculus of Variations, PDE, image processing and fracture mechanics. In some sense 
the theory of sets of finite perimeter is an analogue for sets of what the theory of Sobolev functions is for 
functions. A very useful notion for Sobolev functions is the notion of \em precise representative \rm for a 
function $f\in W^{1,p}$, this is a function $\tilde{f}\in W^{1,p}$ with $\|f-\ti{f}\|_{W^{1,p}}=0$ and 
$\tilde{f}$ has additional smoothness and regularity properties. Another very useful result is Whitney's theorem 
that for any $\ep>0$ there exists a smooth function $\hat{f}$ such that $\|\hat{f}-f\|_{W^{1,p}}<\ep$. We prove 
analogues results for sets of finite perimeter in the plane and we apply them to 
give a short proof that quasiminimizing sets are $BV_l$ extension domains. 

Recall a set $E$ is said to be of \em finite perimeter \rm in domain $\Omega$ if it is measurable and  
$$
\mathrm{Per}(E,\Omega):=\sup\lt\{\int_{E} div \phi:\phi\in \lt[C^1_c(\Omega)\rt]^N, \|\phi\|_{\infty}\leq 1\rt\}<\infty.
$$
A set is simply called a \em set of finite perimeter \rm  if it is of finite perimeter in $\R^n$ and we define $\mathrm{Per}(E):=\mathrm{Per}(E,\R^n)$. When dealing with sets of finite perimeter it is common to work with the 
representative that has negligible points removed. Given a set $E$ we define $\mathrm{Neg}(E):=\lt\{x\in E: \lt|B_r(x)\cap E\rt|=0\text{ for some }r>0\rt\}$. Let $\wt{E}=E\backslash \mathrm{Neg}(E)$. By Proposition 3.3 \cite{amb} we know 
$\|\chara_{E}-\chara_{\wt{E}}\|_{BV(\R^2)}=0$. If a set $F$ is such that $\mathrm{Neg}(F)=\emptyset$ 
we say $F$ is \em shaved. \rm 

A set of finite perimeter is called \em indecomposable \rm iff for any disjoint subsets $A, B\subset E$ such that $E=A\cup B$ 
we have $\mathrm{Per}(E)=\mathrm{Per}(A)+\mathrm{Per}(B)$ then either $\lt|A\rt|=0$ or $\lt|B\rt|=0$. 

The main theorem we will establish in this note is the following. 

\begin{a2}
\label{TT1} 
If $F\subset \R^2$ is a shaved indecomposable set then the closure $\bar{F}$ of $F$ has the property 
$$
\|\chara_{F}-\chara_{\bar{F}}\|_{BV(\R^2)}=0.
$$ 
Consequently the characteristic function of any indecomposable set is $BV$ equivalent to the characteristic function 
of a closed set. 
\end{a2}

A straightforward corollary to this is: 
\begin{a6}
\label{TT} 
Suppose $E\subset \R^2$ is a set of finite perimeter, then for any $\ep>0$ we can find a closed subset $\Eb\subset E$ with finitely many indecomposable components such that  $H^1(\partial^M \Eb\backslash \partial^M E)=0$ and 
\begin{equation}
\label{vic1}
\|\chara_E-\chara_{\Eb}\|_{BV(\R^2)}<\ep.
\end{equation}
\end{a6}

We will show by example that Theorem \ref{TT1} is false for indecomposable sets of finite perimeter in $\R^n$ for 
$n\geq 3$. Our example also shows that Theorem 7 of \cite{amb5} for saturated indecomposable sets in the plane 
does not hold true in dimension three and above. It is unclear to us if Corollary \ref{TT} is true in dimension 
three and above.

\begin{a2}
\label{thmexample}
There exists a shaved set of finite perimeter $S\subset Q=(-1,1)\times (-1,1)\times (-1,1)$ with the following properties
\begin{enumerate}[(i)]
\item \label{ts1} $S$ is connected and $S^c$ is connected. Hence $S$ is an indecomposable saturated set. 
\item \label{ts2} $\lt|\bar{S}\backslash S\rt|>\frac{3999}{4000}$. 
\item \label{ts3} $H^{2}\lt(\partial^M  S\backslash \phi(S^2)\rt)>0$ for any Lipschitz map $\phi:S^2\rightarrow \R^3$.  
\end{enumerate}
\end{a2}

Our main application of Theorem \ref{TT1} will be to show that local quasiminimizers in the plane are $BV_l$ extension domains. 
This is already known as a consequence of the work of David-Semmes, \cite{david} however our proof is much shorter. 
Specifically we say a set $E$ of finite perimeter $E$ is a $K$-\em quasiminimal \rm set in $\Omega$ iff for all 
open $U\subset \subset \Omega$ and all Borel sets $F,G\subset U$ we have 
$\mathrm{Per}(E,U)\leq K \mathrm{Per}((E\cup F)\backslash G,U)$. 
And we say a set $E$ of finite perimeter is locally $K$-quasiminimizing if there exists $\delta>0$ such that 
for any $x\in \partial E$ the set $E$ is a $K$-quasiminimal in $B_{\delta}(x)$. In the case where $E$ is bounded 
$\delta$ can be chosen depending on $x$.

Finally a set $E$ of finite perimeter is a $BV_l$ extension domain if and only if there are constants $c\geq 1$ and $\delta>0$ such
that whenever $u\in BV(E)$ is such that the diameter of the support of $u$
is smaller than $\delta$, then there is a function $Tu\in BV$ such that $\|D Tu\|\leq c\|Du\|(E)$
and $Tu=u$ on $E$.

\begin{a6}
\label{TT2}
If $E\subset \R^2$ is a locally $K$-quasiminimizing set then it is a $BV_l$ extension domain. 
\end{a6}

Note also that if Corollary \ref{TT} were true in dimension three and above then the proof of 
Corollary \ref{TT2} would work in these dimensions too. A generalization of Corollary \ref{TT} to higher dimensions could 
potentially be a useful technical tool in the study of sets of finite perimeter. 

In addition we will obtain the following corollary which is also an easy corollary to Theorem 7 \cite{amb5}. 
Firstly some definitions, Ambrosio et al.\ \cite{amb5} (Definition 2) define a \em hole \rm of a set of finite perimeter 
$S$ to be an indecomposable component of $S^c$ with finite measure. A set $S$ is called \em saturated \rm 
(again see Definition 2, \cite{amb5}) if it is the union of itself and all its holes.  
\begin{a6}[To Theorem \ref{TT1}]
\label{CC1}
Suppose $S\subset \R^2$ is a indecomposable saturated set then there exists an open set $\Sbb$ 
such that $\|\chara_{S}-\chara_{\Sbb}\|_{BV(\R^2)}=0$.
\end{a6}
\bf Acknowledgements. \rm I would thank N.\ Shanmugalingam; the application to quasiminimizing 
sets was worked out through discussions with her. In addition I thank R.L.\ Jerrard who 
essentially provide me with the example constructed in Theorem \ref{thmexample}. Finally I 
would like to thank the referee for careful reading and numerous excellent suggestions that 
greatly improved the readability of this note.

%
%

\section{Sketch of proof of main theorem} 
\label{sketmain}

The proof of the Theorem \ref{TT1} follows from three basic steps. Each follows from the 
last in a fairly natural way. We will firstly state the steps then sketch the reasons they 
hold afterwards. 

Since $E$ is a shaved indecomposable set it has the property 
\begin{equation}
\label{nik30}
\lt|B_r(x)\cap E\rt|>0\text{ for any }x\in E, r>0. 
\end{equation}

\bf Step 1. \rm Let $Z:=\lt\{z\in \mathbb{R}^2: \limsup_{r\rightarrow 0} 
\frac{\lt|E\cap B_r(x)\rt|}{\pi r^2}\geq \frac{1}{2}\rt\}$. We will show we can find a countable 
collection of balls $\lt\{B_{r_n}(x_n):x_n\in Z\rt\}$ with the following properties. 
\begin{enumerate}[(i)]
\item $\lt\{B_{\frac{r_n}{5}}(x_n):n\in \mathbb{N}\rt\}$ are disjoint. 
\item $Z\subset \bigcup_{n} B_{r_n}(x_n)$
\item $\lt|B_{r_n}(x_n)\cap E\rt|\geq \frac{\pi r_n^2}{4}$ for each $n$.
\item $\Pi:=\bigcup_{n} B_{r_n}(x_n)$ is connected.
\end{enumerate}

\bf Step 2. \rm We will show that for $H^1$ a.e.\ $x\in \partial E\backslash Z$ there exists 
$r_x>0$ such that 
\begin{equation}
\label{nik35}
H^1(\partial E\cap B_r(x))>\frac{r}{1600}\text{ for all }r\in \lt(0,r_x\rt].
\end{equation}

\bf Step 3. \rm We will show that $H^1(\partial E\backslash Z)=0$. \nl

-----------------------

\bf Sketch of Step 1. \rm By definition of $Z$ for any $z\in Z$ we can find 
$r_z>0$ such that $\lt|B_{r}(z)\cap Z\rt|\geq \frac{\pi r^2}{4}$ for all $r\in \lt(0,r_z\rt]$. 
So by the $5r$ Covering Theorem we can find a subcollection 
$\lt\{B_{r_{z_n}}(z_n):n\in \mathbb{N}\rt\}$ such that $\lt\{B_{\frac{r_{z_n}}{5}}(z_n):n\in \mathbb{N}\rt\}$ are pairwise disjoint and $Z\subset \bigcup_{n} B_{r_{z_n}}(z_n)$. The only remaining issue is to show that $\Pi:=\bigcup_{n} B_{r_{z_n}}(z_n)$ is 
connected. Suppose it is not, so there are two non-empty disjoint connect components of 
$\Pi_0$ and $\Pi_1$ such that $\Pi=\Pi_0\cup \Pi_1$. Letting $Q_0=E\cap \Pi_0$ and 
$Q_1=E\cap \Pi_1$ it is possible to show that $H^1(\partial^M Q_0\cap \partial^M Q_1)=0$ 
and this implies $E$ is not an indecomposable set, contradiction. Thus $\Pi$ is 
connected. \nl

\bf Sketch of Step 2. \rm Now we can assume without loss of generality that $x\in E^0$. So 
we can find $p_x>0$ such that $\lt|E\cap B_r(x)\rt|<\frac{r^2}{10000}$ for all 
$r\in (0,p_x)$. Note also by (\ref{nik30}) for any $x\in \partial E$ we have that 
$\lt|E\cap B_{\frac{r}{1000}}(x)\rt|>0$ for any $r>0$. So let $r\in (0,p_x)$, by 
Step 1 we have a countable collection of balls $\lt\{B_{r_n}(x_n):n\in \mathbb{N}\rt\}$ 
that satisfy $(i),(ii),(iii),(iv)$. A subcollection of these balls 
$\lt\{B_{r_{p_k}}(x_{p_k}):k\in \mathbb{N}\rt\}$ is such that $B_{r_{p_k}}(x_{p_k})\cap 
B_r(x)\not=\emptyset$ for any $k\in \mathbb{N}$, $E\cap B_r(x)\subset 
\bigcup_{k} B_{r_{p_k}}(x_{p_k})$ and $\wt{\Pi}:=\bigcup_{k} B_{r_{p_k}}(x_{p_k})$ is a 
connected set. So $\wt{\Pi}$ is a `tentacle' of balls that reaches from the outside 
of $B_r(x)$ to $B_{\frac{r}{1000}}(x)$. And note that any ball $B_{r_{p_k}}(x_{p_k})$ has 
at least a quarter of its area is filled by $E$. On the other hand most of $B_r(x)$ is 
empty of $E$. So assume for the moment for simplicity that the tentacle reaches into 
$B_{\frac{r}{1000}}(x)$ in something like a line, pick a direction $v$ that is roughly orthogonal to the line. Now there must be a large set of lines in direction $v$ running through $B_r(x)$ that start at some point in $E\cap B_r(x)$ and end in some point 
in $E^c\cap B_r(x)$. The variation of $\chara_E$ restricted to any of these lines is at 
least $1$, so integrating across these lines gives (\ref{nik35}). Note the reason 
our proof works in $\R^2$ and does not work in higher dimension\footnote{And recall indeed by the example 
constructed in Theorem \ref{thmexample} shows the result in false in higher dimension.} 
is that a `tentacle' in higher dimension has arbitrarily small surface area, where 
as in two dimension the surface area of a tentacle is $O(1)$.  \nl

\bf Sketch of Step 3. \rm 
Suppose $H^1(\partial E\backslash Z)>0$. Assume for simplicity $H^1(\partial E\backslash Z)<\infty$. Let 
$\mu(A):=H^1(\partial^M E\cup (\partial E\backslash Z))$. We can find an open set $U$ with 
$\partial^M E\subset U$ and a compact set $C\subset \partial^M E$ such that 
$\mu(U)\leq \mu(\partial^M E)+\ep$ and $\mu(C)\geq \mu(\partial^M E)-\ep$. We can assume $\ep$ is 
sufficiently small so that $\mu(\partial E\backslash U)>\frac{H^1(\partial E\backslash Z)}{2}$. 

Now since $\delta=\mathrm{dist}(\partial U,C)>0$ we can find a countable collection of 
pairwise disjoint balls $\lt\{B_{r_n}(x_n):x_n\in \partial E\backslash U, r_n<\delta \rt\}$ such that 
$\sum_{n} r_n\geq c_0 \mu(\partial E\backslash U)$ some constant $c_0>0$. By Step 2 we have that 
$\sum_{n} H^1(\partial^M E\cap B_{r_n}(x_n))\geq \frac{c_0}{1600}\mu(\partial E\backslash U)$. Since 
$B_{r_n}(x_n)\cap C=\emptyset$ for all $n$ we have 
$$
\sum_{n} H^1(\partial^M E\cap B_{r_n}(x_n))+H^1(C)\geq \frac{c_0}{1600}\mu(\partial E\backslash U)
+H^1(\partial^M E)-\ep
$$
which is contradiction for small enough $\ep$.

%
%

\section{Sketch of proof of the application to quasiminimizing sets} 
\label{sket}

As stated the main application of Theorem \ref{TT1} is Corollary \ref{TT2}. So to establish this 
we will use the criteria for $BV_l$ extension domains found in \cite{baldi}, \cite{burago}. Namely $E$ is a 
$BV_l$ extension domain 
if for every set of finite perimeter $F\subset E$ with $\mathrm{diam}(F)<\gamma$ we can find $\hat{F}$ with $F\subset \hat{F}$ 
and $\mathrm{Per}(\hat{F},\R^2)\leq C\mathrm{Per}(F,E)$. Following the method of (the preprint form of) 
\cite{ln} we take $\hat{F}$ to be equal to $F$ and it will 
suffice to show that 
\begin{equation}
\label{jjk1}
\mathrm{Per}(\hat{F},\R^2)\leq (1+K)\mathrm{Per}(F,E)
\end{equation}
for any set of finite perimeter $F\subset E$ with $\mathrm{diam}(F)<\gamma$. We will achieve this in 
the following way: For any $\delta>0$ we will find an open set $\Omega$ with $F\subset \subset \Omega$ and 
\begin{equation}
\label{jjk2}
H^1(\partial^M E\cap \Omega)\leq H^1(\partial^M E\cap \partial^M F)+\delta.
\end{equation}
Then by the fact that $E$ is a local $K$-quasiminimizer $\mathrm{Per}(E,\Omega)\leq K\mathrm{Per}(E\backslash F,\Omega)$. 
Now 
\begin{eqnarray}
\label{jpk1}
H^1(\partial^M E\cap \partial^M F)&\leq& H^1(\partial^M E\cap \Omega)\nn\\
&=&\mathrm{Per}(E,\Omega)\nn\\
&\leq&K \mathrm{Per}(E\backslash F,\Omega)\nn\\
&=&K\lt(H^1(\partial^M F\cap E)+H^1\lt(\lt(\partial^M E\backslash \partial^M F\rt)\cap \Omega\rt) \rt)\nn\\
&\overset{(\ref{jjk2})}{\leq}& K \mathrm{Per}(F,E)+K\delta.
\end{eqnarray}
As this holds for arbitrary $\delta>0$ we have inequality (\ref{jjk1}). Note this inequality can not work unless 
we can show (\ref{jjk2}) and hence establish $H^1((\partial^M E\backslash \partial^M F)\cap \Omega)<\delta$. Now 
for arbitrary sets of finite perimeter $F\subset E$ it is not true we can find $\Omega$ such that $F\subset \subset \Omega$ 
and (\ref{jjk2}) holds true. For a counter example let $E=B_1(0)$ and pick $z_0\in \partial E$, let $\alpha\in (0,1)$ and let 
$\zeta_k$ be the set of points with rational coordinates in $E\cap B_{\alpha}(z_0)$. Then the set 
$F=\bigcup_{k} B_{\alpha 2^{-(1000+n)}}(\zeta_n)\cap E$ is a set of finite perimeter for which (\ref{jjk2}) is false for any open set $\Omega$ with $F\subset \subset \Omega$. 

However we will be able to carry out this argument by replacing $F$ by the set $\Fb$ afforded 
to us by Corollary \ref{TT}. The set $\Fb$ has almost the same characteristics of $F$ and we can 
find an open set $\Omega$ with $\Fb\subset \subset \Omega$ such that $H^1(\partial^M E\cap \Omega)\leq H^1(\partial^M E\cap \partial^M \Fb)+\delta$. 
So we can carry out the chain of inequalities to establish (\ref{jpk1}). Having established 
(\ref{jpk1}) for $\Fb$ the same inequality follows for $F$ with arbitrarily small error by (\ref{vic1}). See Lemma 
\ref{LA1} and Section \ref{S5.2} for full details.

\section{Preliminaries}

As in Definition 3.60 \cite{amb} we every $t\in \lt[0,1\rt]$ we let $E^t$ denote the set of points of $t$ density, i.e. 
\begin{equation}
\label{neq2}
E^t:=\lt\{x\in \R^n: \lim_{r\rightarrow 0} \frac{\lt|E\cap B_{r}(x)\rt|}{\Gamma(n) r^n}=t\rt\}
\end{equation} 
where $\Gamma(n)=\lt|B_1(0)\rt|$.

It is a fundamental result of Federer (see Theorem 3.61 \cite{amb}) that if $E$ is a set of finite perimeter in $\Omega$, 
$H^{n-1}\lt(\Omega\backslash \lt(E^0\cup E^{\frac{1}{2}}\cup E^1\rt)\rt)=0$. 

The measure theoretic boundary is defined by 
$$
\partial^M E:=\lt\{x\in \R^n:\limsup_{r\rightarrow 0} \frac{\lt|E\cap B_{r}(x)\rt|}{\lt|B_r(x)\rt|}>0\text{ and } \limsup_{r\rightarrow 0} \frac{\lt|E^c\cap B_{r}(x)\rt|}{\lt|B_r(x)\rt|}>0\rt\}.
$$
It is well known that $H^{n-1}(\partial^M E\backslash E^{\frac{1}{2}})=0$ and 
$\lt|D\chara_E\rt|=H^{n-1}_{\partial^M E}$, see Theorem 3.59, Theorem 3.61 \cite{amb}.

\section{Preliminary lemmas}

\begin{a1}
\label{LL2}
Suppose $E\subset \R^2$ is indecomposable. Then we can find a countable collection 
$$
\lt\{B_{r_n}(x_n):x_n\in E\rt\}
$$ 
such that 
\begin{equation}
\label{nik50}
\lt(\bigcup_{n} B_{r_n}(x_n)\rt)\text{ is connected},
\end{equation}
\begin{equation}
\label{nik51}
\lt|B_{r}(x_n)\cap E\rt|>\frac{\pi r^2}{4}\text{ for all }r\in \lt(0,r_{n}\rt], n\in \mathbb{N}.
\end{equation}

As a consequence for $H^1$ a.e.\ $x\in \partial E\backslash (E^{1}\cup \partial^M E)$ there exists $r_x>0$ such that 
\begin{equation}
\label{qeq1}
H^1(\partial E\cap B_{r}(x))>\frac{r}{1600}\text{ for all }r\in  \lt(0,r_{x}\rt].
\end{equation}

\end{a1}

\em Proof of Lemma \ref{LL2}. \rm Let $Z=E^{1}\cup E^{\frac{1}{2}}$. By Theorem 3.61 \cite{amb}
\begin{equation}
\label{neq2.5}
H^1(E\backslash Z)=0.
\end{equation}

So for any $x\in Z$ there exists $r_x>0$ such that 
\begin{equation}
\label{sseq1}
\lt|B_{r}(x)\cap E\rt|>\frac{\pi r^2}{4}\text{ for all }r\in \lt(0,r_{x}\rt]. 
\end{equation}

\em Step 1. \rm By the $5r$ Covering Theorem (see Theorem 2.11 \cite{mat}) we can find a disjoint sub-collection $\lt\{B_{\frac{r_{x_n}}{5}}(x_n):x_n\in Z\rt\}$ such that $Z\subset \bigcup_{n} B_{r_{x_n}}(x_n)=:\Pi$. We will show $\Pi$ is connected. 

\em Proof of Step 1. \rm We argue by contradiction. Suppose $\Pi$ is disconnected. Let $\Pi_0$ be a connected non-empty component of $\Pi$ and let $\Pi_1=\Pi\backslash \Pi_0$. Now 
define $Q_0=E\cap \Pi_0$ and $Q_1=E\cap \Pi_1$. These are both the intersection of two sets of finite perimeter 
and hence are sets of finite perimeter. 

We claim 
\begin{equation}
\label{qeq2}
H^1\lt(\lt(E\backslash Q_0\rt)\backslash Q_1\rt)=0.
\end{equation}
So let $\omega=E\backslash Z$, note $H^1(\omega)\overset{(\ref{neq2.5})}{=}0$ and $E\backslash (Q_0\cup Q_1)=E\backslash \Pi\subset (Z\cup \omega)\backslash \Pi\subset \omega$. So equality (\ref{qeq2}) follows. We also claim 
\begin{equation}
\label{qeq3}
Q_1\backslash (E\backslash Q_0)=\emptyset.
\end{equation}
Now $(Q_1\cup Q_0)\backslash E=\Pi\cap E\backslash E=\emptyset$ so (\ref{qeq3}) is immediate. 
So we actually have 
\begin{equation}
\label{qeq4}
H^1\lt(Q_1\triangle (E\backslash Q_0)\rt)=0.
\end{equation}
By Proposition 3.38 \cite{amb} this is more than enough to conclude 
\begin{equation}
\label{nik88}
\mathrm{Per}(Q_1)=\mathrm{Per}(E\backslash Q_0).
\end{equation}
We will show 
\begin{equation}
\label{qeq5}
H^1\lt(Q_0^{\frac{1}{2}}\cap Q_1^{\frac{1}{2}}\rt)=0.
\end{equation}
Now if $x\in Q_1^{\frac{1}{2}}$ we must have 
$$
\lim_{r\rightarrow 0} \frac{\lt|E\cap \Pi_1\cap B_r(x)\rt|}{\pi r^2}=\frac{1}{2}
$$ 
so 
\begin{equation}
\label{nik1.5}
\lim_{r\rightarrow 0} \frac{\lt|E\cap B_r(x)\rt|}{\pi r^2}\geq \frac{1}{2}.
\end{equation}
Hence (\ref{nik1.5}) together with Theorem 3.61 \cite{amb} implies 
$$
H^1\lt(Q_1^{\frac{1}{2}}\backslash Z\rt)=0
$$
so 
\begin{equation}
\label{sseq1.6}
H^1(Q_1^{\frac{1}{2}}\backslash \Pi)=0. 
\end{equation}
Now note that if 
$x\in Q_1^{\frac{1}{2}}\cap \Pi$ we can not have $x\in \Pi_0$ because $\Pi_0$ is open and so $B_{\delta}(x)\subset \Pi_0$ 
for some $\delta>0$. By the fact $x\in Q_1^{\frac{1}{2}}$ also we must be able to find $y\in B_{\delta}(x)\cap Q_1$ 
which contradicts the fact $\Pi_1$, $\Pi_0$ are disjoint. Thus 
\begin{equation}
\label{nik2}
Q_1^{\frac{1}{2}}\cap \Pi\subset Q_1^{\frac{1}{2}}\cap \Pi_1.
\end{equation}  

And in the same way since $\Pi_1$ is open for all small enough $r$ we have that $B_r(x)\subset \Pi_1$. 
Thus if $x\in Q_1^{\frac{1}{2}}\cap \Pi_1$ by definition of $Q_1^{\frac{1}{2}}$, 
$\lim_{r\rightarrow 0} \frac{\lt|(E\cap \Pi_1)^c\cap B_r(x)\rt|}{\pi r^2}=\frac{1}{2}$ 
we actually have 
$$
\lim_{r\rightarrow 0} \frac{\lt|E^c\cap B_r(x)\rt|}{\pi r^2}=\frac{1}{2} 
$$ 
hence $x\in E^{\frac{1}{2}}$. Thus
\begin{equation}
\label{nik3}
Q_1^{\frac{1}{2}}\cap \Pi_1\subset E^{\frac{1}{2}}.
\end{equation}
Hence
\begin{eqnarray}
\label{fzf1}
&~&
H^1(Q_1^{\frac{1}{2}}\backslash (E^{\frac{1}{2}}\cap \Pi_1))\nn\\
&~&\qd\qd \leq 
H^1(Q_1^{\frac{1}{2}}\backslash \Pi)+H^1(Q_1^{\frac{1}{2}}\cap \Pi\backslash Q_1^{\frac{1}{2}}\cap \Pi_1)
+H^1(Q_1^{\frac{1}{2}}\cap \Pi_1\backslash E^{\frac{1}{2}}\cap \Pi_1)\nn\\
&~&\qd\qd \overset{(\ref{sseq1.6}),(\ref{nik2}),(\ref{nik3})}{=}0.
\end{eqnarray}
Now going in the opposite direction if $x\in E^{\frac{1}{2}}\cap \Pi_1$ again 
since $\Pi_1$ is open for all small enough $r$ we have that 
$$
\lim_{r\rightarrow 0} \frac{\lt|B_r(x)\cap E\rt|}{\pi r^2}=\lim_{r\rightarrow 0} \frac{\lt|B_r(x)\cap E\cap \Pi_1\rt|}{\pi r^2}=\frac{1}{2}
$$
and 
$$
\lim_{r\rightarrow 0} \frac{\lt|B_r(x)\cap E^c\rt|}{\pi r^2}=\lim_{r\rightarrow 0} \frac{\lt|B_r(x)\cap (E\cap \Pi_1)^c\rt|}{\pi r^2}=\frac{1}{2}.
$$
So $E^{\frac{1}{2}}\cap \Pi_1\subset (E\cap \Pi_1)^{\frac{1}{2}}=Q_1^{\frac{1}{2}}$ hence putting this together with 
(\ref{fzf1}) we have established 
$$
H^1\lt(E^{\frac{1}{2}}\cap \Pi_1\triangle Q_1^{\frac{1}{2}}\rt)=0. 
$$
In exactly the same way we can show that $H^1(E^{\frac{1}{2}}\cap \Pi_2\triangle Q_2^{\frac{1}{2}})=0$. Since $\Pi_1$, $\Pi_0$ are disjoint 
this completes the proof of (\ref{qeq5}). 

Now $E=Q_0\cup Q_1$ and $\mathrm{Per}(Q_0)=H^1(Q_0^{\frac{1}{2}})$, 
$\mathrm{Per}(Q_1)=H^1(Q_1^{\frac{1}{2}})$ 
and $\mathrm{Per}(E)=\mathrm{Per}(Q_0)+\mathrm{Per}(Q_1)\overset{(\ref{nik88})}{=}\mathrm{Per}(Q_0)+\mathrm{Per}(E\backslash Q_0)$. So as $Q_0$, $E\backslash Q_0$ are both 
sets of finite perimeter this contradicts the fact $E$ is indecomposable. This concludes the proof of Step 1. \nl

\em Step 2. \rm We will establish (\ref{qeq1}).

\em Proof of Step 2. \rm Firstly by Theorem 3.61 \cite{amb} we can assume $x\in E^0$. 
Now recall from the sketch of the proof (see property (\ref{nik30})), for any $y\in E$ we have $\lt|B_r(y)\cap E\rt|>0$ for any $r>0$. Thus since $x\in \partial E$ it has the same property. Now since $x\in E^0$ we can find $p_x>0$ such that 
\begin{equation}
\label{fkke7}
\lt|E\cap B_r(x)\rt|<\frac{r^2}{100000}\text{ for all }r\in (0,p_x). 
\end{equation}
However we must also have that $\lt|E\cap B_{\frac{r}{1000}}(x)\rt|>0$. So 
by property (i) we have established in Step 1 we can find a countable collection 
$$
\lt\{B_{r_k}(x_k):x_k\in E, r_k<\frac{r}{1000}\rt\} 
$$
such that $E^{1}\cup E^{\frac{1}{2}}\subset \bigcup_{k} B_{r_k}(x_k)$ and $\lt(\bigcup_{k} B_{r_k}(x_k)\rt)$ is connected. Now 
pick a point $x_0\in E^{1}\cap B_{\frac{r}{1000}}(x)$ and a point $y_0\in E^{1}\backslash B_r(x)$. Since $\bigcup_{k} B_{r_k}(x_k)$ 
is open and connected it is path connected and so 
we must be able to find a path $\phi:\lt[0,t\rt]\rightarrow \bigcup_{k} B_{r_k}(x_k)$ with $\phi(0)=x_0$, $\phi(t)=y_0$. 
Let $s\in (0,t)$ be the smallest number such that $\phi(s)\in \partial B_{\frac{r}{4}}(x)$, by compactness 
clearly this number exists. So 
\begin{equation}
\label{eqaak1}
\lt\{\phi(w):0\leq w\leq s\rt\}\subset B_{\frac{r}{4}}(x).
\end{equation}
Let $v=\frac{\phi(s)-x_0}{\lt|\phi(s)-x_0\rt|}$. For any vector $w$ let $\la w \ra :=\lt\{\lm w: \lm \in \R\rt\}$. 
Define $P_{V}(x)$ to be the orthogonal projection of $x$ onto subspace $V$. Note that for every $k$
\begin{equation}
\label{wzw1}
H^1\lt(P_{\la v\ra}\lt(B_{r}(x_{k})\cap E\rt)\rt)\geq \frac{r}{16}\text{ for all }r\in \lt(0,r_k\rt]
\end{equation}
since if this was not true we 
would have that $\lt|B_{r}(x_{k})\cap E\rt|\leq \frac{r^2}{8}$ which contradicts (\ref{nik51}).

We know $\phi(\lt[0,t\rt])\subset \bigcup_{k} B_{r_k}(x_k)$, let $\lt\{B_{r_{p_k}}(x_{p_k}):k\in \mathbb{N}\rt\}$ 
be a subcollection defined by 
\begin{equation}
\label{nik301}
\lt\{\phi(w):0\leq w\leq s\rt\}\cap B_{r_{p_k}}(x_{p_k})\not=\emptyset\text{ for any }k.
\end{equation}
So 
\begin{equation}
\label{nik5}
\lt\{\phi(w):0\leq w\leq s\rt\}\subset \bigcup_{k} B_{r_{p_k}}(x_{p_k}).
\end{equation}
And 
\begin{equation}
\label{eqaak2}
\bigcup_{k} B_{r_{p_k}}(x_{p_k})\overset{(\ref{nik301}),(\ref{eqaak1})}{\subset} B_{\frac{r}{3}}(x).
\end{equation}
Now since 
\begin{equation}
\label{sseq3}
P_{\la v\ra}\lt(\lt\{\phi(w):w\in \lt[0,s\rt]\rt\}\rt)\overset{(\ref{nik5})}{\subset} \bigcup_{k} P_{\la v\ra}\lt(B_{r_{p_k}}(x_{p_k})\rt)
\end{equation}
by the $5r$ Covering Theorem (see Theorem 2.11 \cite{mat}) 
we can find a subcollection 
$$
\lt\{B_{r_{q_k}}(x_{q_k}):k\in \mathbb{N}\rt\}
$$ 
such that   
\begin{equation}
\label{sseq4}
\lt\{P_{\la v\ra}\lt(B_{r_{p_k}}(x_{p_k})\rt):k\in \mathbb{N}\rt\}\subset \lt\{P_{\la v\ra}\lt(B_{r_{q_k}}(x_{q_k})\rt):k\in \mathbb{N}\rt\}
\end{equation}
and 
\begin{equation}
\label{sseq2}
\lt\{P_{\la v\ra}\lt(B_{\frac{r_{q_k}}{5}}(x_{q_k})\rt):k\in \mathbb{N}\rt\}\text{ are disjoint}.
\end{equation}

Now as  
$H^1\lt(P_{\la v \ra}\lt(\lt\{\phi(w):w\in \lt[0,s\rt]\rt\}\rt)\rt)\geq \frac{r}{5}$ so putting this together with 
(\ref{sseq3}), (\ref{sseq2}) we have 
\begin{eqnarray}
\label{sseq5}
\sum_{k} r_{q_k}&\overset{(\ref{sseq2})}{=}&\sum_{k} 2^{-1}H^1\lt(P_{\la v\ra}(B_{r_{p_k}}(x_{p_k}))\rt)\nn\\
&\overset{(\ref{sseq3})}{\geq}&2^{-1}H^1\lt(P_{\la v\ra}(\lt\{\phi(w):w\in \lt[0,s\rt]\rt\})\rt)\nn\\
&\geq& \frac{r}{10}. 
\end{eqnarray}
So let $\OI:=\bigcup_{k} B_{\frac{r_{q_k}}{5}}(x_{q_k})\cap E$. By (\ref{wzw1}), (\ref{sseq2}) and (\ref{sseq5}) we have 
\begin{eqnarray}
\label{nik81}
H^1\lt(P_{\la v\ra}(\OI)\rt)&\overset{(\ref{sseq2})}{=}&\sum_{k} H^1\lt(P_{\la v\ra}\lt(B_{\frac{r_{q_k}}{5}}(x_{q_k})\cap E\rt)\rt)\nn\\
&\overset{(\ref{wzw1})}{\geq}&\sum_{k} \frac{r_{q_k}}{80}\nn\\
&\overset{(\ref{sseq5})}{\geq}& \frac{r}{800}.  
\end{eqnarray}
Now by (\ref{eqaak2}) 
$\OI\subset B_{\frac{r}{3}}(x)$. We claim we can find a subset $\Pi\subset P_{\la v \ra}(\OI)$ with 
$\lt|\Pi\rt|\geq \frac{r}{1600}$ such that 
\begin{equation}
\label{fkke8}
H^1\lt(P_{\la v\ra}^{-1}(\omega)\cap E^c\cap A\lt(x,\frac{r}{3},r\rt)\rt)>0\text{ for any }\omega\in \Pi. 
\end{equation}
Suppose this is not true. So there is a set $\Lambda\subset P_{\la v\ra}(\OI)$ where $\lt|\Lambda\rt|\geq \frac{r}{1600}$ 
such that $H^1(P^{-1}_{\la v\ra}(\omega)\cap E^c\cap A(0,\frac{r}{3},r))=0$ for all $\omega\in \Lambda$. 
Let $\Xi=\bigcup_{\omega\in \Lambda} P^{-1}_{\la v\ra}(\omega)\cap E\cap A(0,\frac{r}{3},r)$ and by Fubini 
$\lt|\Xi\rt|\geq \lt|\Lambda\rt|\frac{2r}{3}\geq \frac{r^2}{2400}$ and $\lt|\Xi\cap E^c\rt|=0$ and this 
contradicts (\ref{fkke7}). We have established (\ref{fkke8}). 

Now by Theorem 3.103 \cite{amb} we know that 
\begin{eqnarray}
\frac{r}{1600}&\leq& \int_{w\in \Pi} V\lt(\chara_{E},P_{\la v\ra}^{-1}(w)\cap B_r(x)\rt) dw\nn\\
&\leq& V\lt(\chara_E,B_r(x)\rt)\nn\\
&=&H^1\lt(\partial^M E,B_r(x)\rt).\nn
\end{eqnarray}
So this establishes (\ref{qeq1}). $\Box$\nl

\em Proof of Theorem \ref{TT1}. \rm Firstly as before, without loss of generality 
we can assume that for any $x\in F$, $\lt|F\cap B_{\delta}(x)\rt|>0$ for all $\delta>0$.\nl

\em Step 1. \rm We will show $H^1(\bar{F}\backslash (\partial^MF\cup F^1))=0$.

\em Proof of Step 1. \rm Suppose 
\begin{equation}
H^1(\bar{F}\backslash (\partial^MF\cup F^1))>0.
\end{equation}
Let $Z=\bar{F}\backslash (\partial^M F\cup F^1)$. Note that $Z\cap \mathrm{Int}(F)=\emptyset$ because if 
$x\in Z\cap  \mathrm{Int}(F)$ then $x\in F^1$ which contradicts the definition of $Z$. So as $Z\subset \bar{F}$ we know that $Z\subset \partial F$ and thus 
$Z\subset \partial F\backslash \partial^M  F$. 

Let $\ep>0$. If $H^1(Z)=\infty$ 
pick $B\subset Z$ with $H^{1}(B)=1$ and define 
$S=B\cup \partial^M  F$ otherwise define $S=(\bar{F}\backslash F^1)\cup \partial^M  F$. Let 
$\mu(A):=H^{1}(A\cap S)$. Note 
\begin{equation}
\label{eqq2}
\mu(\partial F\backslash \partial^M  F)\geq 
\min\lt\{H^1(\bar{F}\backslash (F^1\cup \partial^M  F)),H^1(B)\rt\}=\beta>0.
\end{equation}

Measure $\mu$ is Radon so we can find an open set 
$U$ such that $\partial^M  F\subset U$ such that 
\begin{equation}
\label{eqakk8}
\mu(U)<\mu(\partial^M F)+\ep 
\end{equation}
and we can find a compact set $C\subset \partial^M  F$ such that 
\begin{equation}
\label{jo2}
\mu(C)>\mu(\partial^M  F)-\ep.
\end{equation}
And so 
\begin{equation}
\label{jo1}
\mathrm{dist}(C,\partial U)=\delta>0. 
\end{equation}
We can take $\varpi>0$ and a subset $\Gamma_0\subset \partial F\backslash U$ with 
\begin{equation}
\label{eqakk6}
H^1((\partial F\backslash U)\backslash \Gamma_0)<\ep
\end{equation}
 and for any $r\in (0,\varpi), x\in \Gamma_0$ we have 
that $\mu(B_r(x))=H^1(B_r(x)\cap S)\leq 2r$. By Lemma \ref{LL2} we can find $\sigma>0$, $\Gamma_1\subset \Gamma_0$ such that 
\begin{equation}
\label{nik1}
H^1(\partial^M F\cap B_r(x))>\frac{r}{1600}\text{ for all }x\in \Gamma_1, r\in (0,\sigma).
\end{equation}
And 
\begin{equation}
\label{eqakk7}
H^1(\Gamma_0\backslash \Gamma_1)<\ep. 
\end{equation}

Now by Vitali covering theorem (see Theorem 2.8 \cite{mat}) we can find 
a pairwise disjoint collection 
$$
\lt\{B_{r_k}(x_k):x_k\in  \Gamma_1 \rt\}
$$
such that 
\begin{equation}
\label{eqakk5}
\mu(\Gamma_1\backslash (\bigcup_{k=1}^{\infty} B_{r_k}(x_k)))=0 
\end{equation}
and 
\begin{equation}
\label{equak22}
\mathrm{sup}\lt\{r_k:k\in \N\rt\}<\min\lt\{\frac{\delta}{2},\varpi,\sigma\rt\}. 
\end{equation}
Note 
\begin{eqnarray}
\beta&\leq&\mu(\partial F\backslash \partial^M F)\nn\\
&=&\mu(\partial F\backslash U)+
\mu(\partial F\cap U\backslash \partial^M F)\nn\\
&\leq& \mu(\partial F\backslash U)+\mu(U\backslash \partial^M F)\nn\\
&\overset{(\ref{eqakk8})}{\leq}&\mu(\partial F\backslash U)+\ep.\nn
\end{eqnarray}
So $\mu(\partial F\backslash U)>\beta-\ep$. Now 
\begin{equation}
\label{eqv2.2}
\sum_{k} \mu(B_{r_k}(x_k))
\overset{(\ref{eqakk5}),(\ref{eqakk7}),(\ref{eqakk6})}{\geq} \mu(\partial F\backslash U)-2\ep\overset{(\ref{eqakk8})}{\geq} \mu(\partial F\backslash \partial^M F)-3\ep.
\end{equation}
By (\ref{nik1}) and the fact $x_k\in \Gamma_1$, $r_k<\sigma$ we have 
$$
1600 H^1(\partial^M  F\cap B_{r_k}(x_k) )\geq r_k\text{ for all }k.
$$ 
Since we choose $r_k<\varpi$ and $x_k\in \Gamma_1\subset \Gamma_0$ be definition of $\Gamma_0$, 
$\mu(B_{r_k}(x_k))\leq 2r_k$ so 
\begin{equation}
\label{eqq1}
\mu(B_{r_k}(x_k))\leq 3200 H^1(\partial^M  F\cap B_{r_k}(x_k))
\end{equation}
thus putting this together with (\ref{eqv2.2}) we have 
\begin{equation}
\label{eqv2}
\sum_{k} 3200 H^1(\partial^M  F\cap B_{r_k}(x_k))  \geq \mu(\partial F\backslash \partial^M F)-3\ep
\overset{(\ref{eqq2})}{\geq}\frac{\beta}{2}.
\end{equation}

Now since $x_k\in \Gamma_1\subset \partial F\backslash U$ by (\ref{jo1}) and (\ref{equak22}) 
we know that 
$C\cap B_{r_k}(x_k)=\emptyset$ for any $k$, so 
\begin{eqnarray}
\label{eqv1}
H^{1}(\partial^M  F)&\geq&H^1(C)+\sum_{k} H^{1}(B_{r_k}(x_k)\cap \partial^M F)\nn\\
&\overset{(\ref{eqv2}),(\ref{jo2})}{\geq}& H^1(\partial^M  F)+\frac{\beta}{6400}-\ep 
\end{eqnarray}
which is a contradiction assuming $\ep$ is small enough. \nl

\em Step 2. \rm We will show that $\bar{F}$ is a set of finite perimeter and $D\chara_{\bar{F}}=D\chara_{F}$. 

\em Proof of Step 2. \rm By Step 1 $\lt|\bar{F}\backslash \partial^M F\cup F^1\rt|=\lt|\bar{F}\backslash F\rt|=0$. So by Proposition 3.38 \cite{amb} we have $\mathrm{Per}(\bar{F})=\mathrm{Per}(F)$ and 
hence $\bar{F}$ is a set of finite perimeter and $D \chara_{\bar{F}}=D \chara_{F}$. $\Box$

--------------------------------------------------------------------------

\section{The Applications}

\subsection{Quasiminimizing sets}

The following lemmas hold true in $\R^n$ without additional complexity, so we state them in $\R^n$.

\begin{a1}
\label{L1}
Given as set of finite perimeter $S$, suppose $H^{n-1}(A)=0$ then $\mathrm{Per}(S,A)=0$.
\end{a1}
\em Proof Lemma \ref{L1}. \rm By Theorem 1.9 (2), Corollary 1.11 \cite{mat} measure $\mu$ defined by 
$$
\mu(H):=H^{n-1}(\partial^M  S\cap H)
$$
is a Radon measure. 

Suppose set $A$ has the property $H^{n-1}(A)=0$. Then $\mu(A)=0$, so 
\begin{eqnarray}
0&=&\inf\lt\{\mu(V):A\subset V, V \text{ is open }\rt\}\nn\\
&=&\inf\lt\{\mathrm{Per}(S,V):A\subset V, V \text{ is open }\rt\}\nn\\
&=&\mathrm{Per}(S,A).\nn
\end{eqnarray}
Using the fact $A\rightarrow \mathrm{Per}(S,A)$ is also Radon measure, see Proposition 3.38(a), Proposition 1.43 \cite{amb}.

%
%

\begin{a1}
\label{L2}
Let $E,F$ be sets of finite perimeter in $\R^n$, $F\subset E$. Then 
\begin{equation}
\label{eqd30}
\mathrm{Per}(F, \partial^M E)=H^{n-1}(\partial^M  F\cap \partial^M E).
\end{equation}
and 
\begin{equation}
\label{eqd10}
\mathrm{Per}(E,\partial^M  F)\geq H^{n-1}(\partial^M  F\cap \partial^M E).
\end{equation}
Hence 
\begin{equation}
\label{eh2}
\mathrm{Per}(F,\partial^M E)\leq\mathrm{Per}(E,\partial^M F).
\end{equation}
\end{a1}
 
\em Proof of Lemma \ref{L2}. \rm Note that the measure $\nu(A):=H^{n-1}(\partial^M  F\cap A)$ is a Radon measure. Letting $\sigma>0$ be some small number. 
Pick open $\Omega$ with 
\begin{equation}
\label{nik23}
\partial^M E\cap \partial^M F\subset \Omega 
\end{equation}
such that 
\begin{equation}
\label{eqb2}
\nu(\Omega)=H^{n-1}(\partial^M  F\cap \Omega)<H^{n-1}(\partial^M E\cap \partial^M  F)+\sigma=\nu(\partial^M E)+\sigma. 
\end{equation}

\em Step 1. \rm We will establish (\ref{eqd30}).

\em Proof of Step 1. \rm Define 
\begin{equation}
\label{eqa20}
B:=\lt\{x\in (\partial^M  F)^c: \limsup_{r\rightarrow 0}\frac{H^{n-1}(\partial^M  F\cap B_r(x))}{r^{n-1}}>0  \rt\}
\end{equation}
and 
$$
D:=\lt\{x\in \partial^M E: \liminf_{r\rightarrow 0} H^{n-1}(\partial^M E\cap B_r(x))/r^{n-1}<1\rt\}.
$$ 
By Theorem 6.2 \cite{mat} we have that $H^{n-1}(B)=0$ and since $\partial^M E$ is rectifiable by 
Theorem 16.2 \cite{mat} $H^{n-1}(D)=0$. 

So using Lemma \ref{L1} for the last equality 
\begin{eqnarray}
\label{eqb7}
\mathrm{Per}(F, \partial^M E)&=&\mathrm{Per}(F, \partial^M E\cap \partial^M  F)+
\mathrm{Per}(F, \partial^M E\backslash \partial^M  F)\nn\\
&=&\mathrm{Per}(F, \partial^M E\cap \partial^M  F)+
\mathrm{Per}(F, \partial^M E\backslash (\partial^M  F\cup B\cup D)).
\end{eqnarray}  

Note $H^{n-1}(\partial^M E)<C$. Let $\ep>0$. We can find a decreasing sequence of number $\delta_m\rightarrow 0$ such that the sets 
\begin{equation}
\label{nik11}
\UI_m:=\lt\{\begin{array}{lcl} x\in \partial^M E\backslash (\partial^M F\cup B\cup D)&:&\frac{H^{n-1}(\partial^M  F\cap B_r(x))}{r^{n-1}}<\ep\\
&~&\text{ and }\frac{H^{n-1}(\partial^M E\cap B_r(x))}{r^{n-1}}\geq \frac{1}{2} \text{ for all }r\in (0,\delta_m)
\end{array} \rt\}
\end{equation}
have the property that $\partial^M E\backslash (\partial^M F\cup B\cup D)\subset \bigcup_{m=1}^{\infty} \UI_m$. 

Let $U_1=\UI_1\backslash \lt(\bigcup_{i=2}^{\infty} \UI_i\rt)$,  
$U_2=\UI_2\backslash \lt(\bigcup_{i=3}^{\infty} \UI_i\rt)$, ... $U_k=\UI_k\backslash \lt(\bigcup_{i=k+1}^{\infty} \UI_i\rt)$. 
Now $U_l$, $U_k$ are disjoint for any $l,k$ and 
\begin{equation}
\label{nik60}
\partial^M E\backslash (\partial^M F\cup B\cup D)=\bigcup_{i=1}^{\infty} \UI_i=
\bigcup_{i=1}^{\infty} U_i.
\end{equation} 
Since $\lt\{U_1,U_2,\dots \rt\}$ are pairwise disjoint 
\begin{equation}
\label{jj12}
\sum_{m} H^{n-1}(U_m)\overset{(\ref{nik60})}{<}H^{n-1}(\partial^M E\backslash \partial^M F).
\end{equation}

Pick $m$. Since by Section 5.1 \cite{mat} we have $S^1(U_m)\leq 2 H^1(U_m)$  where $S^1$ denotes 
Spherical Hausdorff measure. So we can find a collection $\lt\{B_{\frac{r_k}{2}}(z_k):r_k<\delta_m \rt\}$ such that 
$U_m\subset \bigcup_{k} B_{\frac{r_k}{2}}(z_k)$ and $\sum_{k} \Gamma(n-1)\lt(\frac{r_k}{2}\rt)^{n-1}\leq 2 S^1(U_m)$. Now 
for each $k$ we can pick $x_k\in B_{\frac{r_k}{2}}(z_k)\cap U_m$ and then we have a collection of balls 
\begin{equation}
\label{nik101}
\lt\{B_{r_k}(x_k):x_k\in U_m, r_k<\delta_m\rt\}
\end{equation}
such that 
\begin{equation}
\label{nik61}
U_m\subset \bigcup_{k} B_{r_k}(x_k)=:V_m
\end{equation}
and 
\begin{equation}
\label{nik12}
\sum_{k} r_k^{n-1}\leq c S^1(U_m)\leq c H^1(U_m).
\end{equation}
Hence using the fact that $x_k\in U_m\subset \UI_m$ we have 
\begin{eqnarray}
\label{wj1}
H^{n-1}(V_m\cap \partial^M  F)&\overset{(\ref{nik61})}{\leq}&\sum_{k} H^{n-1}(\partial^M  F \cap B_{r_k}(x_k))\nn\\
&\overset{(\ref{nik101}),(\ref{nik11})}{\leq}& c\sum_{k} \ep r_{k}^{n-1}\nn\\
&\overset{(\ref{nik12})}{\leq}& c\ep H^{n-1}(U_m).
\end{eqnarray}
Now 
\begin{eqnarray}
\mathrm{Per}(F,\partial^M E\backslash (\partial^M F\cup B\cup D))&\overset{(\ref{nik60}),(\ref{nik61})}{\leq}& \mathrm{Per}\lt(F,\bigcup_{m} V_m\rt)\nn\\
&\leq& \sum_{m} H^{n-1}(V_m\cap \partial^M  F)\nn\\
&\overset{(\ref{wj1})}{\leq}&  \sum_{m} c\ep H^{n-1}(U_m)\nn\\
&\overset{(\ref{jj12})}{\leq}& c\ep H^{n-1}(\partial^M E\backslash \partial^M F). 
\end{eqnarray}
Taking the limit as $\ep\rightarrow 0$ we have $\mathrm{Per}(F,\partial^M E\backslash (\partial^M F\cup B\cup D))=0$. So 
putting this together (\ref{eqb7}) we have 
\begin{equation}
\label{ea50}
\mathrm{Per}(F, \partial^M E)=\mathrm{Per}(F, \partial^M E\cap \partial^M  F).
\end{equation}
And
\begin{eqnarray}
\label{ea}
\mathrm{Per}(F,\partial^M E\cap \partial^M  F)&\overset{(\ref{nik23})}{\leq}& \mathrm{Per}(F,\Omega)\nn\\
&\leq& H^{n-1}(\partial^M  F\cap \Omega)\nn\\
&\overset{(\ref{eqb2})}{\leq}&  H^{n-1}(\partial^M  F\cap \partial^M E)+\sigma.
\end{eqnarray}
Now as $\sigma$ is arbitrary, from (\ref{ea}) and (\ref{ea50}), $\mathrm{Per}(F,\partial^M E)\leq 
H^{n-1}(\partial^M F\cap \partial^M E)$. Conversely for any open set $\Omega$ with $\partial^M E\subset \Omega$ 
we have $\mathrm{Per}(F,\Omega)\geq H^{n-1}(\partial^M F\cap \partial^M E)$. Thus by taking the 
infimum over all such open sets we have $\mathrm{Per}(F,\partial^M E)\geq H^{n-1}(\partial^M F\cap \partial^M E)$ 
and this completes the proof of Step 1. \nl

\em Step 2. \rm We will establish (\ref{eqd10}).

\em Proof of Step 2. \rm First note that $\partial^M E\cap \partial^M  F$ is 
an $(n-1)$ rectifiable set. Let $\delta>0$, for $H^{n-1}$ a.e.\ $x\in \partial^M E\cap \partial^M  F$ there exists $r_x>0$ such that 
\begin{equation}
\label{eqd1}
\lt|H^{n-1}(\partial^M F\cap \partial^M E\cap B_h(x))-n\Gamma(n)h^{n-1}\rt|<\delta h^{n-1}\text{ for any }h\in (0,r_x). 
\end{equation}
Let $\alpha>0$. Defining $\Gamma_0:=\lt\{x\in \partial^M F\cap \partial^M E:\alpha<r_x \rt\}$ 
for all small enough $\alpha>0$
we have that 
\begin{equation}
\label{eqd85}
H^{n-1}(\partial^M F\cap \partial^M E\backslash \Gamma_0)<\delta. 
\end{equation}

Since $A\rightarrow  \mathrm{Per}(E,A)$ is a Radon measure we can extract 
a compact subset $\Gamma_1\subset \Gamma_0$ such that 
$$
\mathrm{Per}(E,\Gamma_1)+\delta>\mathrm{Per}(E,\Gamma_0).
$$

And we can find a compact subset $\Gamma_2\subset \Gamma_0$ such that 
$H^{n-1}(\Gamma_2)+\delta>H^{n-1}(\Gamma_0)$. Let $\Gamma_3=\Gamma_1\cup \Gamma_2$. So $\Gamma_3$ is 
compact and has the properties 
\begin{equation}
\label{eqd81}
\mathrm{Per}(E,\Gamma_3)+\delta>\mathrm{Per}(E,\Gamma_0)\text{ and }
H^{n-1}(\Gamma_3)+\delta>H^{n-1}(\Gamma_0).
\end{equation}
And 
\begin{equation}
\label{jk4}
\Gamma_3\subset \Gamma_0\subset \partial^M F\cap \partial^M E.
\end{equation}

Again since we are dealing with a Radon measure we can find an open set $U$ with $\Gamma_3\subset U$ and 
\begin{equation}
\label{edq60}
\mathrm{Per}(E,U)<\mathrm{Per}(E,\Gamma_3)+\delta.
\end{equation}

Now $\mathrm{dist}(\Gamma_3,U)=:\beta>0$ by Vitali covering theorem we can 
find a collection of pairwise disjoint balls 
$\lt\{B_{r_k}(x_k):x_k\in \Gamma_3, k\in \N\rt\}$ where 
$$
\sup\lt\{r_k:k\in \N\rt\}<\min\lt\{\beta,\alpha\rt\}
$$
and 
\begin{equation}
\label{eqd80}
H^{n-1}(\Gamma_3\backslash (\bigcup_{k=1}^{\infty} B_{r_k}(x_k)))=0. 
\end{equation}
Now since $\bigcup_{k=1}^{\infty} B_{r_k}(x_k)\subset U$ we know 
\begin{equation}
\label{eqd4}
\mathrm{Per}\lt(E,\bigcup_{k=1}^{\infty} B_{r_k}(x_k)\rt)
\leq \mathrm{Per}(E,U)\overset{(\ref{edq60})}{<}\mathrm{Per}(E,\Gamma_3)+\delta.
\end{equation}
But as $\Gamma_3\subset \Gamma_0$
\begin{eqnarray}
\label{eqd3}
\mathrm{Per}\lt(E,\bigcup_{k=1}^{\infty} B_{r_k}(x_k)\rt)&=&
\sum_{k=1}^{\infty} \mathrm{Per}\lt(E,B_{r_k}(x_k)\rt)\nn\\
&=&\sum_{k=1}^{\infty} H^{n-1}(\partial^M E\cap B_{r_k}(x_k))\nn\\
&\overset{(\ref{eqd1})}{\geq}&\sum_{k=1}^{\infty} (n\Gamma(n)-\delta)r_k^{n-1}.
\end{eqnarray}
Now recall $\Gamma_3\subset \Gamma_0\subset \partial^M E\cap \partial^M  F$ so  
\begin{eqnarray}
\label{eqd5}
\sum_{k=1}^{\infty} n\Gamma(n)(1+\delta)r_k^{n-1}&\overset{(\ref{eqd1})}{\geq}& 
\sum_{k=1}^{\infty} H^{n-1}(\partial^M E\cap \partial^M  F\cap B_{r_k}(x_k))\nn\\
&\overset{(\ref{jk4})}{\geq}& \sum_{k=1}^{\infty} H^{n-1}(\Gamma_3\cap B_{r_k}(x_k))\nn\\
&\overset{(\ref{eqd80})}{=}&H^{n-1}(\Gamma_3)\nn\\
&\overset{(\ref{eqd81})}{\geq}& H^{n-1}(\Gamma_0)-\delta\nn\\
& \overset{(\ref{eqd85})}{\geq} & H^{n-1}(\partial^M E\cap \partial^M  F)-2\delta.
\end{eqnarray} 
Now putting (\ref{eqd4}), (\ref{eqd3}), (\ref{eqd5}) together we have 
\begin{equation}
\label{eqd6}
\mathrm{Per}\lt(E,\Gamma_3\rt)\geq H^{n-1}(\partial^M E\cap \partial^M  F)-3\delta
-c\sum_{k=1}^{\infty} \delta r_k^{n-1}.
\end{equation}

Now since the collection of balls $\lt\{B_{r_k}(x_k):x_k\in \Gamma_3, k\in \N\rt\}$ are pairwise disjoint and 
$x_k\in \Gamma_3$ we know by (\ref{eqd1}) we know $\mathrm{Per}(E,B_{r_k}(x_k))\geq \frac{n\Gamma(n) r_k^{n-1}}{2}$
\begin{eqnarray}
\frac{n\Gamma(n)}{2} \sum_{k} r_k^{n-1}&\leq& \sum_{k} \mathrm{Per}(E,B_{r_k}(x_k))\nn\\
&\leq&\mathrm{Per}(E,\R^n).\nn
\end{eqnarray}
So 
\begin{equation}
\label{eqd2}
\sum_{k} r_k^{n-1}\leq c\mathrm{Per}(E,\R^n).
\end{equation}
Now putting (\ref{eqd6}) and (\ref{eqd2}) together we have 
\begin{equation}
\label{eqd9}
\mathrm{Per}\lt(E,\Gamma_3\rt)\geq H^{n-1}(\partial^M E\cap \partial^M  F)-3\delta
-c\delta\mathrm{Per}(E,\R^n).
\end{equation}
Since $\mathrm{Per}\lt(E, \partial^M  F\rt)\geq\mathrm{Per}\lt(E,\Gamma_3\rt)$ and 
$\delta$ is arbitrarily small this establishes (\ref{eqd10}). \nl 

\em Proof of Lemma completed. \rm By applying Step 1 and Step 2 we 
have 
\begin{eqnarray}
 \mathrm{Per}(E,\partial^M F)&\overset{(\ref{eqd10})}{\geq}&H^{n-1}(\partial^M E\cap \partial^M  F)\nn\\
&\overset{(\ref{eqd30})}{=}&\mathrm{Per}(F,\partial^M E)\;\;\;\; \Box
\end{eqnarray}

%
%

\begin{a1}
\label{L2.5}
We will show that if $E\subset \R^n$ is an open set of finite perimeter and $\FF\subset E$ is a relatively closed set of 
finite perimeter, then   
\begin{equation}
\label{er1}
\partial^M  (E\backslash \FF)=(\partial^M  \FF\cap E)\cup (\partial^M  E\backslash \partial^M  \FF).
\end{equation}
 \end{a1}

\em Proof. \rm The lemma above holds true for arbitrary sets of finite perimeter $E$ and $\FF\subset E$, we will 
only need it for open set $E$ and relatively closed set $\FF$ and as the proof is easier in this case we argue only this 
result. 

Let 
\begin{equation}
\label{er12}
\Gamma_1=\partial^M  (E\backslash \FF)\cap \FF\text{ and }\Gamma_2=\partial^M  (E\backslash \FF)\backslash \FF.
\end{equation}

Since $\FF$ is relatively closed, if $x\in \Gamma_2$ then for some small enough $\delta>0$, 
$(B_{\delta}(x)\cap E)\cap \FF=\emptyset$ so 
$$
\lim_{r\rightarrow 0} \frac{\lt|B_r(x)\cap (E\backslash \FF)\rt|}{\Gamma(n)r^n}=
\lim_{r\rightarrow 0} \frac{\lt|B_r(x)\cap E\rt|}{\Gamma(n)r^n}=\frac{1}{2}
$$
and thus $\Gamma_2\subset \partial^M  E$ and hence $\Gamma_2\subset (\partial^M  E\backslash \FF)\subset 
\partial^M  E\backslash \partial^M  \FF$. On the other hand $\Gamma_1\subset \FF\subset E$ and $E$ is open. By 
definition of $\partial^M (E\backslash \FF)$ for any $x\in \Gamma_1$, for all small enough $r$ we have 
$\lt|B_r(x)\cap (E\backslash \FF)\rt|\approx \frac{\Gamma(n)r^n}{2}$ and $\lt|B_r(x)\cap (E\backslash \FF)^c\rt|\approx \frac{\Gamma(n)r^n}{2}$. But as $E$ is open $B_r(x)\subset E$ for all small enough $r$ so we must have $x\in \partial^M \FF$. Thus $\Gamma_1\subset \partial^M  \FF\cap E$ and hence as $\Gamma_1\cup \Gamma_2
=\partial^M (E\backslash \Fb)$ this establishes (\ref{er1}). $\Box$ \nl

%
%

\begin{a1}
\label{LA1}

Let $E\subset \R^n$ be a locally $K$-quasiminimizing set and let $\FF\subset E$ be a relatively closed subset of finite perimeter. Then 
\begin{equation}
\label{eh60}
\mathrm{Per}(\FF,\partial^M E)\leq K \mathrm{Per}(\FF,E).
\end{equation}
\end{a1}
\em Proof of Lemma \ref{LA1}. \rm First note by Theorem 4.2 \cite{ln} we know the topological boundary $\partial E$ is equal to the measure theoretic 
boundary $\partial^M E$.\nl

\em Step 1. \rm We will show there exists open set $\Omega$ with $\FF\subset \Omega$ such that 
\begin{equation}
\label{rtr1}
H^{n-1}((\partial E\backslash \partial^M  \FF)\cap \Omega)\leq \ep. 
\end{equation}

\em Proof of Step 1. \rm Let $\mu(A):=H^{n-1}(A\cap \partial E)$. So $\mu$ is a Radon measure on $\R^n$ and hence 
we must be able to find open set $U$ with $\partial E\cap \partial \FF\subset U$ such that 
$\mu(U)\leq H^{n-1}(\partial E\cap \partial \FF)+\ep$. 

Note $\partial \Fb\cap \partial E$ is a closed set. So $\delta=\inf\lt\{\lt|x-y\rt|:x\in \partial \FF\cap \partial E, y\not\in U\rt\}>0$. Let $\Pi=\bigcup_{x\in \partial \Fb\cap \partial E} B_{\delta}(x)$. Note 
\begin{equation}
\label{equaz1}
\mu(\Pi)\leq \mu(U)\leq H^{n-1}(\partial E\cap \partial^M F)+\ep. 
\end{equation}
Note by 
compactness there 
exists $\lm>0$ such that $N_{\lm}(\partial E)\cap \Fb\subset \Pi$. Let 
$\Lambda=\bigcup_{x\in \Fb\backslash N_{\lm}(\partial E)} B_{\frac{\lm}{2}}(x)$. Note $\Lambda\cap \partial E=\emptyset$. 
And 
$$
\Fb=(\Fb\cap N_{\lm}(\partial E))\cup (\Fb\backslash N_{\lm}(\partial E))\subset \Pi\cup \Lambda.
$$ 
Now 
\begin{eqnarray}
\mu(\Lambda\cup \Pi)&=&\mu((\Lambda\cup \Pi)\cap \partial E)=\mu(\Pi\cap \partial E)\nn\\
&\overset{(\ref{equaz1})}{\leq}& H^{n-1}(\partial E\cap \partial^M \Fb)+\ep.\nn
\end{eqnarray}
So letting $\Omega=\Lambda\cup \Pi$, open set $\Omega$ satisfies (\ref{rtr1}).\nl

\em Proof of Lemma completed. \rm Note that for any open set $\OI$ with $E\subset \OI$, 
$$
H^{n-1}(\partial^M  \FF\cap E\cap \Omega)\leq H^{n-1}(\partial^M  \FF\cap E)\leq \mathrm{Per}(\FF,\OI).
$$
Thus 
\begin{equation}
\label{eqab10}
\mathrm{Per}(\FF,E)\geq H^{n-1}(\partial^M  \FF\cap E\cap \Omega). 
\end{equation}
So using Lemma \ref{L2}
\begin{eqnarray}
\label{eh1}
\mathrm{Per}(E\backslash \FF,\Omega)&\overset{(\ref{er1})}{=}&H^{1}(\partial^M  \FF\cap E\cap \Omega)+H^{1}((\partial^M E\backslash \partial^M  \FF)\cap \Omega)\nn\\
&\overset{(\ref{eqab10}),(\ref{rtr1})}{\leq}&\mathrm{Per}(\FF,E)+\ep.
\end{eqnarray}
Now since $\FF\subset \subset \Omega$ and $E$ is a quasiminimizer 
\begin{eqnarray}
\label{eh4}
\mathrm{Per}(E,\Omega)&\leq&K\mathrm{Per}(E\backslash \FF,\Omega)\nn\\
&\overset{(\ref{eh1})}{=}&K(\mathrm{Per}(\FF,E)+\ep)
\end{eqnarray}
And since $\FF$ is closed, $\partial^M \FF\subset \FF$ so 
\begin{equation}
\label{nikk1}
\mathrm{Per}(\FF,\partial^M E)\overset{(\ref{eh2})}{\leq} \mathrm{Per}(E,\partial^M \FF)\leq \mathrm{Per}(E,\FF),  
\end{equation}
thus
\begin{eqnarray}
\label{ew1}
\mathrm{Per}(\FF,\partial^M E)+\mathrm{Per}(E,\Omega\backslash \FF)&\overset{(\ref{nikk1})}{\leq}& 
\mathrm{Per}(E,\FF)+\mathrm{Per}(E,\Omega\backslash \FF)\nn\\
&=&\mathrm{Per}(E,\Omega)\nn\\
&\overset{(\ref{eh4})}{\leq}&K \mathrm{Per}(\FF,E)+K\ep.
\end{eqnarray}
Since this gives 
$$
\mathrm{Per}(\FF,\partial^M E)\leq K\mathrm{Per}(\FF,E)+K\ep. 
$$
And as $\ep$ is arbitrary this establishes (\ref{eh60}). $\Box$\nl

%
%
%
%

\subsection{Proof of Theorem \ref{TT2}.} \label{S5.2}

As stated in Section \ref{sket} by the criteria for $BV_l$ extension domain of \cite{baldi}, 
\cite{burago}, $E$ is a $BV_l$ extension domain if for every set of finite perimeter $F\subset E$ with 
$\mathrm{diam}(F)<\delta$ we can find $\hat{F}$ with $F\subset \hat{F}$ with $\mathrm{Per}(\hat{F},\R^n)\leq 
(1+K)\mathrm{Per}(F,E)$. We will take $\hat{F}=F$ and we will show 
\begin{equation}
\label{jw1}
\mathrm{Per}(F,\R^2)\leq (1+K)\mathrm{Per}(F,E).
\end{equation}
So start using Corollary \ref{TT} we can take closed subset $\FF\subset F$ with $H^1(\partial^M \FF\backslash \partial^M F)=0$ and 
\begin{equation}
\label{ddf1}
\|\chara_{\FF}-\chara_{F}\|_{BV}<\ep.
\end{equation}
Now by Lemma \ref{LA1} 
\begin{equation}
\label{ddf2}
\mathrm{Per}(\FF,\partial^M E)\leq K\mathrm{Per}(\FF,E).
\end{equation}
Note since $\FF$ is closed, $\mathrm{Per}(\FF,\FF^c)=0$. And note since $\FF\subset E$ we know 
$H^1(E^0\cap \FF)=0$ and so $\mathrm{Per}(\FF,E^0\cap \FF)=0$ by Lemma \ref{L1}. So 
\begin{eqnarray}
\label{jjww1}
\mathrm{Per}(\FF,E^0)=\mathrm{Per}(\FF,E^0\cap \FF^c)+\mathrm{Per}(\FF,E^0\cap \FF)=0.
\end{eqnarray}
 
Now again by Lemma \ref{L1}
\begin{eqnarray}
\label{jjww2}
\mathrm{Per}(\FF,\R^2)&\leq& \mathrm{Per}(\FF,E^1)+ \mathrm{Per}(\FF,\partial^M E)+ \mathrm{Per}(\FF,E^0)\overset{(\ref{ddf2}),(\ref{jjww1})}{\leq} (1+K)\mathrm{Per}(\FF,E).
\end{eqnarray}
Hence 
\begin{eqnarray}
\label{ddf3}
\mathrm{Per}(F,\R^2)&\overset{(\ref{ddf1})}{\leq}& \mathrm{Per}(\FF,\R^2)+\ep\overset{(\ref{jjww2})}{\leq} (1+K)\mathrm{Per}(\FF,E)+\ep.
\end{eqnarray}
For any open set $\OI$ with $E\subset \OI$ we have 
$$
\mathrm{Per}(\FF,E)\leq \mathrm{Per}(\FF,\OI)=H^1(\partial^M \FF,\OI)\leq H^1(\partial^M F,\OI)=\mathrm{Per}(F,\OI).
$$
As this is true for all open sets $\OI$ with $E\subset \OI$ we have $\mathrm{Per}(\FF,E)\leq \mathrm{Per}(F,E)$ 
thus putting this into (\ref{ddf3}) we have $\mathrm{Per}(F,\R^2)\leq (1+K)\mathrm{Per}(F,E)+\ep$ and 
this holds for all $\ep>0$, so this establishes (\ref{jw1}). $\Box$

%
%
%
%

\subsection{Saturated indecomposable sets.\nl\nl}

\em Proof of Corollary \ref{CC1}. \rm 

\em Step 1. \rm First we show that $\lt|S\rt|<\infty$. 

\em Proof of Step 1. \rm By the isoperimetric inequality (see Theorem 3.46 \cite{amb}) we have that 
$$
\min\lt\{\lt|S\rt|, \lt|\R^2\backslash S\rt|\rt\}\leq (4\pi)^{-1}\lt(\mathrm{Per}(S,\R^2)\rt)^{2} 
$$
Now if $\lt|\R^2\backslash S\rt|\leq (4\pi)^{-1}\lt(\mathrm{Per}(S,\R^2)\rt)^{2}$ then $\R^2\backslash S$ must 
have an indecomposable component with finite measure which contradicts the definition of saturated, so this 
can not happen. \nl

\em Step 2. \rm Let 
$$
a=\inf\lt\{\tau>0:\lt|D\chara_S\rt|(B_{\tau}(0))>0\rt\}\text{ and }b=\sup\lt\{\omega>0:\lt|D\chara_S\rt|(\R^2\backslash B_{\omega}(0))>0\rt\}.
$$
We will show 
\begin{equation}
\label{jkq1}
H^1(\partial^M S\cap \partial B_s)>0\text{ or }H^1(S\cap \partial B_s)>0\text{ for all }s\in (a,b).
\end{equation}

\em Proof of Step 2. \rm Suppose not, so for some $s\in (a,b)$ we have that 
\begin{equation}
\label{edq1}
H^1(\partial^M S\cap \partial B_s)=0\text{ and }H^1(S\cap \partial B_s)=0.
\end{equation}
Now by definition of $a,b$ we know 
$\mathrm{Per}(S,B_s)>0$ and $\mathrm{Per}(S,\bar{B_s}^c)>0$. Now 
let $A=S\cap B_s$ and $B=S\cap \bar{B_s}^c$. So $A$, $B$ are both sets of 
finite perimeter.

Now for $H^1$ a.e.\ $x\in \partial^M A\backslash \partial^M S$, since $A\subset S$ we must have 
$x\in S^1$ (recall definition (\ref{neq2})) and so $x\in \partial B_s$. Thus 
$$
H^1(\partial^M A\backslash (\partial^M S\cap \bar{B_s}\cup (S\cap \partial B_s)))=0.
$$ 
In the same way 
$$
H^1(\partial^M B\backslash (\partial^M S\cap B_s^c\cup (S\cap \partial B_s)))=0.
$$ 
Thus by (\ref{edq1}) we have that $H^1(\partial^M A\cap \partial^M B)=0$ and 
as $\lt|A\cap B\rt|=0$ by Proposition 1 \cite{amb5} we have $\mathrm{Per}(A)+\mathrm{Per}(B)=\mathrm{Per}(S)$ which contradicts the fact that $S$ is indecomposable. 

Now suppose $b=\infty$, then letting 
\begin{equation}
\label{jke1}
H_1:=\lt\{s>a:H^1(\partial^M S\cap \partial B_s)>0\rt\}.
\end{equation}
and  
$$
H_2:=\lt\{s>a:H^1(S\cap \partial B_s)>0\rt\}.
$$

We have two cases to consider. Either $\lt|H_1\rt|=\infty$ or $\lt|H_2\rt|=\infty$. We will deal with each in turn. Firstly 
we will use a standard adaptation of Theorem 3.103 \cite{amb} we have that for any $u\in L^1(B_R(0))$ 
\begin{equation}
\label{llk1}
V(u,B_R(0))=\int_{0}^R V(u,\partial B_s(0)) ds
\end{equation}
where 
\begin{equation}
\label{llk1.6}
V(u,\partial B_s(0))=\sup\lt\{\int_{0}^{2\pi s} u\lt(s\cos\lt(\frac{\theta}{s}\rt),s\sin\lt(\frac{\theta}{s}\rt)\rt)\phi'(\theta) d\theta:\phi\in C_0^{\infty}(\lt[0,2\pi s\rt])\rt\}.
\end{equation}
Now suppose $\lt|H_1\rt|=\infty$ then 
\begin{equation}
\label{nik70}
V(\chara_S,\R^2)\overset{(\ref{llk1}),(\ref{jke1})}{\geq}\int_{H_1} V(\chara_S,\partial B_s(0)) ds=\infty. 
\end{equation}
Which contradicts the fact that $S$ a set of finite perimeter. Now suppose $\lt|H_2\rt|=\infty$. 
Let $B_2:=\lt\{s\in H_2: H^1(\partial B_s(0)\cap S)\geq \pi s\rt\}$ then by the Coarea formula 
\begin{eqnarray}
\infty>\lt|S\rt|\geq \int_{s\in B_2} H^1(\partial B_s(0)\cap S) ds\geq \int_{s\in B_2} \pi s ds\geq \lt|B_2\cap \lt\{x:x\geq 1\rt\}\rt|.
\end{eqnarray}
Hence $\lt|B_2\rt|<\infty$. Now for any $t\in H_2\backslash B_2$ we must have  $V(\chara_S,\partial B_s(0))\geq 1$. 
Thus as in (\ref{nik70}), $V(\chara_S,\R^2)\geq \lt|H_2\backslash B_2\rt|=\infty$ which contradicts the fact $S$ has finite 
perimeter. So in either case we have a contradiction, thus $b<\infty$. \nl

\em Step 3. \rm We will show $\lt|\overline{B_{b}(0)}^c\cap S\rt|=0$. 

\em Proof of Step 3. \rm Since $\lt|D\chara_S\rt|(\overline{B_{b}(0)}^c)=0$ either 
$\lt|\overline{B_{b}(0)}^c \cap S\rt|=0$ or $\lt|\overline{B_{b}(0)}^c \cap S^c\rt|=0$ in the latter case $S^c$ has an indecomposable component of finite measure 
which contradicts the fact that $S$ is stratified. \nl 

\em Step 4. \rm We will show $S^c$ has only one indecomposable component.

\em Proof of Step 4. \rm Let $U_1,U_2,\dots $ be the indecomposable components of $S^c$. Since $S$ is saturated, $\lt|U_k\rt|=\infty$ 
for any $k\in \mathbb{N}$. Thus $\lt|U_k\cap B_b^{c}\rt|>0$ for all $k$. Since $H^1(\partial^M U_k\cap B_b^c)=0$, thus 
$\lt|D\chara_{U_k}\rt|(B_b^c)=0$ so either $\chara_{U_k}\equiv 1$ on $B_b^c$ or $\chara_{U_k}\equiv 0$ on $B_b^c$. As $\lt|U_k\cap B_b^c\rt|>0$ the former is true and so $\lt|B_b^c\backslash U_k\rt|=0$. Thus for every $k,l$ we have 
$\lt|U_k\cap U_l\rt|>0$ which is a contradiction. So there is only one indecomposable component. \nl

\em Proof of Corollary completed. \rm By Theorem \ref{TT1}, $S^c$ can be 
represented by an indecomposable component $\widetilde{\Pi}$ which is closed. So 
$\Sbb:=(\widetilde{\Pi})^c$ then $\Sbb$ is open and $D\chara_S=D\chara_{\Sbb}$. $\Box$

\section{Example: Proof of Theorem \ref{thmexample}}.

\em{The construction}. \rm Let $\lt\{\zeta_k:k\in \mathbb{N}\rt\}$ be an enumeration of the points in the 
unit square in the $e_1,e_2$ plane with rational coordinates, i.e.\ points in the set 
$Q_1:=\lt\{xe_1+ye_2: x,y\in \lt[-1,1\rt]\rt\}$ that can be written as $x e_1+y e_2$ for some $x,y\in \lt[-1,1\rt]\cap \mathbb{Q}$. 

We will define a thin column in direction $e_3$ around $\zeta_k$ by 
$$
\Pi_k:=\lt\{x e_1+y e_2+ z e_3+\zeta_k: \sqrt{(x-\zeta_k\cdot e_1)^2+(y-\zeta_k\cdot e_2)^2}<2^{-10000-k}, 0< z\leq 1\rt\}.
$$
And we define $\OI=\bigcup_{k\in \mathbb{N}} \Pi_k$. Let 
$$
\UUI=\lt\{(x,y,z):-1\leq x\leq 1, -1\leq y \leq 1, -1 \leq z\leq 0 \rt\}. 
$$
Finally let 
$$
\SI=\OI\cup \UUI.
$$
Since $\OI$ and $\UUI$ are shaved so $\UUI$ is also a shaved set.\nl

\em Step 1. \rm We will show $\SI$ is path connected. 

\em Proof of Step 1. \rm From the construction it is clear that 
for any $p\in \SI$, the line interval $\lt[p, (p\cdot e_1) e_1+(p\cdot e_2) e_2-\frac{e_3}{2} \rt]
\subset \UUI$. Thus for any $p_1,p_2\in \OI$
\begin{equation}
\label{uk1}
L_1=\lt[p_1, (p_1\cdot e_1) e_1+(p_1\cdot e_2) e_2-\frac{e_3}{2} \rt]
\subset \SI
\end{equation}
and 
\begin{equation}
\label{uk2}
L_2=\lt[p_2, (p_2\cdot e_1) e_1+(p_1\cdot e_2) e_2-\frac{e_3}{2} \rt]
\subset \SI.
\end{equation}
Finally 
\begin{equation}
\label{uk3}
L_3=\lt[ (p_1\cdot e_1) e_1+(p_1\cdot e_2) e_2-\frac{e_3}{2} , (p_2\cdot e_1) e_1+(p_1\cdot e_2) e_2-\frac{e_3}{2} \rt]
\subset \UUI\subset \SI.
\end{equation}
So 
$$
L_1\cup L_2\cup L_3\subset \SI.
$$
Thus $p_1$ and $p_2$ are connected in $\SI$. \nl

\em Step 2. \rm We will show $\SI^c$ is connected.

\em Proof of Step 2. \rm Let $Q=\lt(-1,1\rt)\times \lt(-1,1\rt)\times \lt(-1,1\rt)$. If 
$p\in Q\backslash \SI$ then $p+\lm e_3\not \in \SI$ for any $\lm>0$. Thus 
given any two points $p_1,p_2\in Q\backslash \SI$ we have 
$\lt[p_1,p_2+2e_3\rt]\subset \SI^c$ and $\lt[p_1,p_2+2e_3\rt]\subset \SI^c$ and 
thus the path 
$$
\PPI=\lt[p_1,p_1+2e_3\rt]\cup \lt[p_1+2e_3,p_2+2e_3\rt]\cup \lt[p_2+2e_3,p_2\rt]\subset \SI^c,
$$
thus $\PPI$ is connected. \nl

\em Step 3. \rm We will show $\SI$ is a set of finite perimeter. 

\em Proof of Step 3. \rm Let $\OI_N=\bigcup_{k=1}^{N} \Pi_k$ and 
define 
$$
\SI_N=\OI_N\cup \UUI.
$$
Note that $\SI_N$ is a set of finite perimeter and 
$$
\mathrm{Per}(\SI_N)\leq \sum_{k=1}^N 2\pi 2^{-10000-k}+17\leq 18.
$$
Now $\SI_N$ converges in measure to $\SI$ (see Remark 3.36 \cite{amb}) so by Proposition 3.37 (b) we 
know $\mathrm{Per}(\SI)\leq 18$ and thus $\SI$ is a set of finite perimeter. \nl

\em Step 4. \rm For a.e.\ $p\in Q\backslash \SI$ there exists $r_p>0$ such that 
\begin{equation}
\label{jb0}
\lt|\SI\cap B_r(p)\rt|>0,\; \lt|\SI^c\cap B_r(p)\rt|>0 \text{ for all }r\in (0,r_p) 
\end{equation}
and 
\begin{equation}
\label{jb1}
\mathrm{Per}(\SI, B_r(p))>0\text{ for any }r\in (0,r_p).
\end{equation}

\em Proof of Step 4. \rm By Lebesgue density theorem for a.e. $p\in Q\backslash \SI$ there exists 
$r_p>0$ such that 
\begin{equation}
\label{jb3}
\frac{\lt|B_r(p)\cap \SI^c\rt|}{\lt|B_r(p)\rt|}>\frac{1}{2}\text{ for any }r\in \lt(0,r_p\rt]. 
\end{equation}
Now for any $r\in \lt(0,r_p\rt]$ we can find $x,y\in \mathbb{Q}$ such that 
$(x,y,p\cdot e_3)\in \OI\cap B_{\frac{r}{4}}(p)$ and as $\OI$ is an open set for some 
$\rho\in (0,\frac{r}{4})$ we know 
\begin{equation}
\label{jb2}
\lt|B_{\rho}(p)\cap \OI\rt|>0. 
\end{equation}
So (\ref{jb2}) and (\ref{jb3}) together establish (\ref{jb0}). 

Now arguing by contradiction and assuming $\mathrm{Per}(\SI\cap B_r(x))=0$ then as 
this means $\lt|D\chara_{\SI}\rt|(B_r(x))=0$ so by Theorem 3.44, \cite{amb} we have that 
$\int_{B_r(x)} \lt|\chara_{\SI}-(\chara_{\SI})_{B_r(x)}\rt| dz=0$. However by (\ref{jb2}), (\ref{jb3}) and this is a contradiction. 
Thus (\ref{jb1}) is established. \nl

\em Step 5. \rm We will show that $Q\backslash \OI\subset \bar{\SI}$.

\em Proof of Step 5. \rm For any $p\in Q\backslash \OI$ either we have 
$p\in \UI$ or $p\in (Q\backslash \UI)\backslash \SI$. 
So using (\ref{jb0}) of Step 4 for the latter case we know that 
for a.e.\ $p\in Q\backslash \OI$, $p\in \bar{\SI}$. So 
there exists a subset $\GI\subset Q\backslash \OI$ with $\lt|(Q\backslash \OI)\backslash \GI\rt|=0$ and 
$\GI\subset \bar{\SI}$. Thus $\bar{\GI}\subset \bar{\SI}$ and since 
$Q\backslash \OI\subset \bar{\GI}$ this completes the proof of Step 5.\nl

\em Proof of Theorem completed. \rm By Step 1 and Step 2, $\SI$ and $\SI^c$ are connected. Since $\SI$ is a set of finite 
perimeter by Proposition 2 \cite{amb5} we have that $\SI$ is an indecomposable set. Since 
$H^2(\partial^M \SI^c)=H^2(\partial^M \SI)<\infty$ in the same way by Proposition 2, $\SI^c$ is indecomposable. Thus it is its own only indecomposable  component. Hence (\ref{ts1}) is established. 

Now note by Step 5, $Q\backslash \OI\subset \bar{\SI}$ and 
$$
\lt|Q\backslash \SI\rt|=4-\sum_{k=1}^{\infty} \pi(2^{-10000-k})^2\geq \frac{3999}{4000},
$$ 
so (\ref{ts2}) is established. 

Finally by Step 4 we can find a set $\GI\subset Q\backslash \SI$ with $\lt|(Q\backslash \SI)\backslash \GI\rt|=0$ 
and for every $p\in \GI$ there exists $r_p>0$ such that (\ref{jb1}) holds true. 
So arguing by contradiction suppose there exists a Lipschitz map $\phi:S^2\rightarrow \R^3$ such that 
$H^2(\partial^M \SI\backslash \phi(S^2))=0$. Now by Theorem 3.59, Theorem 3.61 \cite{amb}, 
$$
H^2(\phi(S^2)\cap B_r(p))\geq H^2(\partial^M \SI\cap B_r(p))
=\mathrm{Per}(\SI,B_r(p))>0\text{ for all }r\in (0,r_p).
$$ 
So $p\in \bar{\phi(S^2)}=\phi(S^2)$. Thus $\GI\subset \phi(S^2)$ but this is a contradiction because 
$H^2(\phi(S^2))\leq (\mathrm{Lip}(\phi))^2 H^2(S^2)= 4\pi(\mathrm{Lip}(\phi))^2$. Thus we have established 
(\ref{ts3}) and completed the proof of the theorem. $\Box$\nl

\em Remark. \rm The example constructed in Theorem \ref{thmexample} also shows that Theorem 7 of \cite{amb5} has no analogue in 
higher dimension. Note by property (\ref{ts1}) there are no bounded components of $\SI^c$ so $\SI$ has 
no holes and is therefore saturated. Thus the set $\SI$ constructed is an example of a set Ambrosio et al.\ define as \em simple \rm (see Definition 3, Definition 2 \cite{amb5}). And by property (\ref{ts3}) $\partial^M S$ can not be almost everywhere 
covered by the Lipschitz image of the sphere.

\end{document}